\documentclass[12pt,twoside]{article}
\usepackage[all]{xy}
\usepackage{amssymb}
\usepackage{modif} 

\font\teneusm=eusm10 at 12pt
\font\seveneusm=eusm7 at 9pt
\font\fiveeusm=eusm5 at 7pt
\newfam\eusmfam
\textfont\eusmfam=\teneusm
\scriptfont\eusmfam=\seveneusm
\scriptscriptfont\eusmfam=\fiveeusm

\font\teneufm=eufm10 at 12pt
\font\seveneufm=eufm7 at 9pt
\font\fiveeufm=eufm5 at 7pt
\newfam\eufmfam
\textfont\eufmfam=\teneufm
\scriptfont\eufmfam=\seveneufm
\scriptscriptfont\eufmfam=\fiveeufm
\def\eufm{\teneufm\fam\eufmfam}

\font\tenmsbm=msbm10 at 12pt
\font\sevenmsbm=msbm7 at 9pt
\font\fivemsbm=msbm5 at 7pt
\newfam\msbmfam
\textfont\msbmfam=\tenmsbm
\scriptfont\msbmfam=\sevenmsbm
\scriptscriptfont\msbmfam=\fivemsbm
\def\msbm{\tenmsbm\fam\msbmfam}

\def\theoremunit{subsubsection}
\newtheorem{theorem}[\theoremunit]{Theorem}
\newtheorem{proposition}[\theoremunit]{Proposition}
\newtheorem{corollary}[\theoremunit]{Corollary}
\newtheorem{lemma}[\theoremunit]{Lemma}
\newtheorem{definition}[\theoremunit]{Definition}

\newtheorem{void}[\theoremunit]{\setbox0=\hbox{\ }\kern-\wd0}

\newenvironment{remark}{\begin{void}{\bf Remark}\rm\ }
{\end{void}}

\def\proofbox{\null\nobreak\hfill\rule{6pt}{6pt}}

\newcounter{acounter}
\newenvironment{enum}%
{\begin{list}{(\roman{acounter})}{\usecounter{acounter}}}%
{\end{list}}
{\begin{list}{\alph{acounter}.}{\usecounter{acounter}}}%
{\end{list}}

\newenvironment{claim}[1]{\begin{trivlist} \item[] {\bf #1}\em }
{\end{trivlist}}

\newenvironment{proof}[1]{\begin{trivlist} \item[] {\bf Proof of~#1\ }}
{\proofbox\end{trivlist}}

\def\hide#1{}

\def\idle{}

\def\onrightpage{\ifodd\value{page}\idle\else
\mbox{}\thispagestyle{empty}\clearpage\fi}

\def\lastleftpage{\mbox{}\ifodd\value{page}
\mbox{}\thispagestyle{empty}\clearpage\fi}
\def\theoremunit{subsection}

\def\Aff{{\msbm A}}

\def\Hom{{\rm Hom}}
\def\sHom{{\underline{\rm Hom}}}
\def\dim#1{{\rm dim}\left(#1\right)}
\def\codim#1{{\rm codim}\left(#1\right)}
\def\H#1#2{{\rm H}_{#1}\!\left(#2\right)}

\def\cH#1#2{{\rm H}^{#1}\!\left(#2\right)}

\def\Gm{{\msbm G}_{m}}

\def\C{{\msbm C}}

\def\field{{\mathrm{\bf k}}}

\def\Spec{{\mathrm{Spec}}}

\def\O{\mathcal{O}}
\def\Ker{{\mathrm{Ker}}}
\def\Coker{{\mathrm{Coker}}}

\def\dd#1{{#1}^{\scriptscriptstyle\vee\!\vee}}
\def\red{{\mathrm{red}}}
\def\sing{{\mathrm{sing}}}
\def\smooth{{\mathrm{smth}}}

\def\Tr{{\mathrm{Tr}}}
\def\GL#1{\mathrm{GL}(#1)}
\def\SL#1{\mathrm{SL}(#1)}

\def\etimes{{\mathbin{\boxtimes}}}
\def\d{{\mathrm d}}

\def\lieg{{\eufm g}}
\def\liet{{\eufm t}}
\def\lieh{{\eufm h}}
\def\quot{\kern-.2em\big/\kern-.3em\big/\kern-.1em}
\def\inv#1#2{{(#2)}^{#1}}
\def\e{{\mathrm{e}}}
\def\K{{\mathrm{K}}}

\def\lie{{\mathrm L}}

\def\under2#1{\underline{\underline{#1}}}

\def\map#1{\stackrel{#1}{\longrightarrow}}

\renewcommand{\thesection}{\arabic{section}}

\begin{document}

%
\title{Differential forms and smoothness\\ of quotients by reductive 
groups}

\author{Guillaume Jamet}

\date{\today}

\maketitle

\begin{abstract}
    Let $\pi : X\map{}Y$ be a good quotient of a smooth variety $X$ 
    by  a reductive algebraic group $G$ and $1\leq k\leq\dim Y$ an integer. We prove 
    that if, locally, any invariant horizontal differential $k$-form on $X$ (resp. any 
    regular differential $k$-form on $Y$)  is a K\"ahler differential
    form on $Y$ then $\codim{Y_{\sing}}>k+1$. We also prove that the 
    dualizing sheaf on $Y$ is the sheaf of invariant horizontal $\dim
    Y$-forms.
\end{abstract}


\pagestyle{myheadings}
\markboth{Guillaume Jamet}{Smoothness of quotients by reductive groups}



\section*{Introduction}

Let $\pi : X\map{}Y$ be a good quotient of a smooth variety $X$ 
by  a reductive algebraic group $G$. How one can bound the dimension 
of the singular locus of $Y$? Since there exists no natural embedding 
of $Y$ in some smooth variety, it seems difficult to describe the 
$n$-th Fitting ideal of the sheaf $\Omega_{Y}^{1}$. J.~Fogarty 
suggests a different approach to this problem by raising in \cite{MR90a:14014} the following 
questions (all schemes are assumed 
to be of finite type over a field of characteristic 0)~:
\begin{claim}{Question} Let $G$ be a finite group acting on a smooth 
variety $X$ and $\pi : X\map{}Y$ the quotient. Is the natural morphism
$$
\Omega^{1}_{Y}\map{}\inv{G}{\Omega^{1}_{X}}
$$
surjective if and only if $Y$ is smooth?
\end{claim}
In that article J. Fogarty verifies that the surjectivity condition is indeed 
necessary. He also proves that, when the 
group $G$ is abelian, this condition is sufficient (\cite[Lemma 5]{MR90a:14014}).

Observe that the module $\inv{G}{\Omega^{1}_{X}}$ is naturally 
isomorphic to $\dd{\Omega^{1}_{Y}}$ and, the variety $Y$ being normal, 
also isomorphic to the module $\omega^{1}_{Y}$ of regular 1-forms (cf. appendix
\ref{rardiff}) and to the module $i_{*}\Omega^{1}_{Y_{\smooth}}$ 
(here $i$ denotes the inclusion $Y_{\smooth}\subset y$). 
It is also easily checked that this problem reduces to 
the case where $X$ is a rational representation of $G$. In particular 
when $G\subset\SL{\C^{2}}$, then $Y={\C^{2}}/G$ is a complete 
intersection and one can give an affirmative answer to the question 
above. However, already in dimension 2 (i.e. $G\subset\GL{\C^{2}}$) 
this question appears to be quite tricky.

Recently M.Brion proved the following result~:
\begin{claim}{Theorem ({\cite[Theorem 1]{MR98k:14067}})} Let $G$ be a
reductive algebraic group acting on 
a smooth affine variety $X$, and let $\pi : X\map{}Y$ be the
quotient. If $Y$ is smooth then the natural morphism
\begin{eqnarray*}
\inv{G}{d\pi} & : &\Omega_{Y}\map{}\inv{G}{\Omega_{X,G}}
\end{eqnarray*}
is an isomorphism.
\end{claim}
Here $\inv{G}{\Omega_{X,G}}$ is the differential graded algebra of {\em invariant 
horizontal differential forms} and $\inv{G}{d\pi}$ is the morphism of 
differential graded algebras induced by the cotangent morphism $d\pi$ (see section \ref{hordiff}).
When $G$ is finite, it is isomorphic to $\inv{G}{\Omega_{X}}$.  
This last theorem clearly suggests to reformulate and investigate Fogarty's question in the more 
general context of quotients by reductive groups. 

The main theorem we prove in this paper is the following, 
thus giving a partial answer to Fogarty's question and also a strong 
converse to Brion's theorem~:
\begin{claim}{Theorem (\ref{smcrit1})} Let $G$ be a reductive
algebraic group acting on a smooth affine variety $X$, with quotient 
map $\pi :X\map{} Y$ and let $k$ be an integer with $1\leq k\leq\dim Y$.
The morphism $\inv{G}{d\pi^{k}}$ is
surjective in codimension $k+1$ if and only if $Y$ is smooth in
codimension $k+1$.    
\end{claim}
We stated these results for affine $G$-schemes, but it is easy 
to see that they generalize immediately to the case of good 
quotients (i.e. affine uniform categorical quotient morphisms $\pi : 
X\map{}Y$, with the terminology of \cite{MR86a:14006}).

In the case of finite abelian groups we also prove~:
\begin{claim}{Theorem (\ref{fgroup3})} Let $G$ be a finite abelian group 
acting on a smooth affine scheme $X$ with quotient $\pi : X\map{} Y$
and let $k$ be an integer with $1\leq k\leq\dim X$. The morphism
$\inv{G}{d\pi^{k}}$ is surjective if and only if $Y$ is smooth.
\end{claim}
This improves the previous result of Fogarty and also shows that, with 
the hypothesis of (\ref{smcrit1}) smoothness in codimension $k+1$ 
doesn't imply that $\inv{G}{d\pi^{k}}$ (or $c_{Y}^{k}$, see below) is surjective.

In order to prove these theorems it is important to 
understand how $\inv{G}{\Omega_{X,G}}$ compares to other sheaves of 
differentials on $Y$, in particular to the sheaves $\tilde{\Omega}_{Y}$ and 
$\omega_{Y}$ (respectively, the sheafs of {\em absolutely regular} and {\em 
regular differential forms}. Cf. appendix \ref{rardiff}).
In his article \cite{MR98k:14067}, M. Brion observed that, as a corollary to 
his theorem and under the additional condition that {\em no invariant 
divisors is 
mapped by $\pi$ onto a closed subscheme of codimension $\geq 2$ in 
$Y$}, there are isomorphisms 
$\inv{G}{\Omega_{X,G}}\simeq\dd{\Omega_{Y}}\simeq\omega_{Y}$. 
This comparison problem is also closely related to the more classical problem of describing the 
dualizing sheaf of a quotient (by a reductive group) variety  as a sheaf of 
invariants. 
It  has been extensively studied by F. Knop in \cite{MR90k:14053}, 
but the expression he obtains for $\omega_{Y}^{n}$ (the canonical 
sheaf if $n=\dim Y$) is again dependent on the existence of the 
preceding ``bad divisors".

Here, using a general machinery of K\"ahler (resp. absolutely 
regular) horizontal differential forms (sections \ref{hordiff} and \ref{arhdiff}) we obtain 
the following comparison statement~:
\begin{claim}{Proposition (\ref{arhcomp})} Let $G$ be a reductive algebraic group, $X$ be
a smooth affine $G$-scheme and $\pi : X\map{}Y$ 
the quotient. There is a sequence of inclusions :
$$
\bar\Omega_{Y}\subseteq\tilde\Omega_{Y}\subseteq\inv{G}{\Omega_{X,G}}
\subseteq\omega_{Y} 
$$
which are equalities on the smooth locus of $Y$. 
\end{claim}
This together with a theorem of Boutot (\cite{MR88a:14005}) leads to the following simple description of the dualizing 
sheaf~:
\begin{claim}{Corollary (\ref{arhdual})} Let $G$ be a reductive algebraic group, $X$ be
a smooth affine $G$-scheme with quotient map $\pi : X\map{}Y$ and let $n=\dim Y$. Then 
the dualizing sheaf on $Y$, $\omega^{n}_{Y}$, is isomorphic to 
$\inv{G}{\Omega^{n}_{X,G}}$.
\end{claim}

Our Proposition \ref{arhcomp} also leads to a more intrinsic version of (\ref{smcrit1})~:
\begin{claim}{Theorem (\ref{smcrit2})} Let $Y$ be the quotient of a smooth
affine variety by a reductive algebraic group and let $k$ be an integer with $1\leq k\leq\dim Y$. 
The fundamental class morphism 
$c_{Y}^{k}$ is surjective in codimension $k+1$ 
if and only if $Y$ is smooth in codimension $k+1$.   
\end{claim}
Note that this result apply in particular when $Y$ is a variety with 
toroidal singularities. Indeed, it is proved in \cite{MR95i:14046} that any toric variety 
can be realized as the good quotient of an open subset of an affine 
space $\Aff^{n}$ by a torus. In fact, for quotient by tori, we 
expect that a statement similar to (\ref{fgroup3}) might hold.

A smoothness criterion much like (\ref{smcrit2}) also holds when 
$Y$ is
locally a complete intersection (\cite{MR42:255} or \cite{mythesis}. Note by the way that quotient 
singularities which are complete intersections are ``exceptionnal" and 
must be singular in codimension 2). 
Even more generally, one may conjecture that for a variety $Y$ with 
reasonnable singularities (see \cite[5.22, p107]{MR90a:14021} in appendix \ref{rardiff}) 
{\em $c_{Y}^{k}$ is surjective in codimension $k$ 
if and only if $Y$ is smooth in codimension $k$}\, (the ``$k+1$'' in 
(\ref{smcrit2}) is clearly a gift of the local quasi-homogeneous 
structure). 

Finally, combined with  results of H. Flenner (\cite{MR89j:14001}, and van Sraten-Steenbrink
\cite{MR87j:32025} in the case of isolated singularities) proposition \ref{arhcomp} implies 
that for $0\leq k<\codim{Y_{\sing}}-1$, we have 
$\tilde\Omega^{k}_{Y}\simeq\inv{G}{\Omega^{k}_{X,G}}\simeq\omega^{k}_{Y}$.
However, the following question (as far as we know) remains open : 
{\em Under the hypotheses of (\ref{arhcomp}) do we have in general isomorphisms 
$\tilde\Omega_{Y}\simeq\inv{G}{\Omega_{X,G}}\simeq\omega_{Y}$ or at 
least
$\inv{G}{\Omega_{X,G}}\simeq\omega_{Y}$?}

%

\paragraph*{Acknowledgements} This work reproduces parts of my Ph.D. 
Thesis, worked out at the Institut de Math\'ematiques de Jussieu.
Many thanks to my Ph.D. advisor, C.~Peskine.

\paragraph*{Notation and conventions}
We work over a fixed field $\field$ of characteristic 0 with algebraic
closure $\bar\field$. 
All the schemes we consider are of
finite type over $\field$. For such a scheme $X$, we denote by $\Omega_{X}$ the
differential graded algebra $\oplus_{k\geq 0}\Omega^{k}_{X/\field}$ of
K\"{a}hler differentials, and
write $\Omega^{k}_{X}$ for $\Omega^{k}_{X/\field}$. 

For $G$ an algebraic group and a $G$-scheme $X$, we denote by
$G$-$\O_{X}$-mod the category of $G$-equivariant
$\O_{X}$-modules.


An affine $\Gm$-scheme $X$ is said to be quasi-conical (this is an ugly
terminology, but, we believe it is consistent with the algebraic
definitions of homogeneous and quasi-homogeneous ideals) if $\O_{X}$ is
generated by homogeneous sections of non-negative weights. We recall that
$X$ is said to be conical when  $\O_{X}$ is
generated by homogeneous sections of  weight 1.

By differential operator, we mean differential operator relative to
$\field$ in the sense of \cite[16.8]{MR39:220}.
 
We denote by $\Gamma$ the decreasing filtration by codimension of
the support~: Let $c$ be an integer. For any $\O_{X}$-module $M$ and
$U\subset X$ an open subset, $\Gamma_{c}M(U)$ is the subgroup of $M(U)$
consisting of the sections having support of codimension $\geq c$ in
$X$. We write $\Gamma_{(c)}$ for $\Gamma_{c}/\Gamma_{c+1}$ and $\bar{M}$
for $\Gamma_{(0)}M$. In particular, when $X$ is integral, $\Gamma_{1}M$
is the submodule of torsion elements and $\bar{M}=\Gamma_{(0)}M$ is $M$
modulo torsion. We recall that this filtration is preserved by
differential operators and in particular by $\O_{X}$-linear morphisms.  
These definitions extend to categories of complexes in
the obvious way. 


By a desingularisation of $X$, we
always mean a desingularisation of $X_{\red}$. We take
(\cite{MR99i:14020}) as a general reference for resolution of
singularities, in particular for the existence of equivariant resolutions.


\section{Horizontal differentials}
\label{hordiff}

Let $G$ be an algebraic group, $\lieg$ its Lie algebra considered as a
 $G$-module via the adjoint representation, and $X$ a
$G$-scheme. We will also consider $G$
as a $G$-scheme by the action of $G$ on itself by inner automorphism. 
We have the following diagram of equivariant maps~:
$$
\xymatrix{
G & G\times X \ar[l]^{p}\ar[d]^{q}\ar@<.5ex>[r]^{\mu} & X \ar@<.5ex>[l]^{s} \\
  & X &
}
$$
where $p$ and $q$ are the projections, $\mu$ is the action map and
$s$ is the section of $\mu$ defined by $x\mapsto (e,x)$. This induces
the following diagram of $G$-equivariant coherent modules on $G\times X$~:
$$
\xymatrix{
\mu^{*}\Omega^{1}_{X}\ar[r]^<<<{d\mu}\ar[rd] &
 \Omega^{1}_{G\times X}=p^{*}\Omega^{1}_{G}\oplus
 q^{*}\Omega^{1}_{X}\ar[d]\\
 & p^{*}\Omega^{1}_{G}
}
$$
Taking the pull-back by $s$ of the diagonal morphism above, we obtain a morphism
\begin{eqnarray*}
 d\mu^{1}_{X,G} & : & \Omega^{1}_{X}\map{}s^{*}p^{*}\Omega^{1}_{G}=\lieg^{\vee}\otimes\O_{X}
\end{eqnarray*}
We then define a morphism 
$d\mu_{X,G} : \Omega_{X}\map{}\Omega_{X}\otimes\lieg^{\vee}$ as 
follows
\begin{eqnarray*}
 d\mu^{k}_{X,G} & : &
\Omega^{k}_{X}\map{}\Omega^{k-1}_{X}\otimes\lieg^{\vee} \\
 d\mu^{k}_{X,G}(df_{1}\wedge\ldots\wedge
df_{k}) &= & \sum_{i=1}^{k}(-1)^{k-i}df_{1}\wedge\ldots\wedge\widehat{df_{i}}\wedge\ldots\wedge
df_{k}\otimes d\mu^{1}_{X,G}(df_{i}),
\end{eqnarray*}

\begin{remark} For an alternative and more rigourous construction of the 
morphisms above, using `multilinear homological algebra', we refer 
to \cite{mythesis}.
\end{remark}

\begin{definition} The $G$-equivariant module
$\Omega^{k}_{X,G}=\Ker(d\mu^{k}_{X,G})$ is called the module of
horizontal $k$-forms. We denote by $\Omega_{X,G}$ the graded algebra
$\oplus_{k\geq 0}\Omega^{k}_{X,G}$.
\end{definition}
The sections of $\Omega_{X,G}$ consists of those forms whose interior product with any
vector field induced by the group action vanishes.

The preceeding construction is natural in $X$. Thus, for any equivariant map $f
: X\map{}Y$ the cotangent morphism induces morphisms
$f^{*}\Omega^{k}_{Y,G}\map{}\Omega^{k}_{X,G}$.
It is also clear from
the construction that if the action of $G$ is  trivial then
$d\mu^{1}_{X,G}=0$ and consequently we have
$\Omega^{k}_{X,G}=\Omega^{k}_{X}$. From these remarks, we deduce~:

\begin{proposition}\label{invmorph} Let $\pi : X\map{}Y$ be a
$G$-invariant morphism, then the cotangent morphism  $d\pi : \pi^{*}\Omega_{Y}\map{}\Omega_{X}$ factors through
$\Omega_{X,G}\subset\Omega_{X}$.
\end{proposition}
\begin{remark}\label{torsionkernel} This last proposition applies in
particular when $\pi$ is a categorical 
quotient of $X$. Assume that $X$ is affine and that $G$ is a 
reductive linear group. 
Let $\pi : X\map{}Y$ be the quotient of $X$. By (\ref{invmorph}) there
is a morphism 
$\pi^{*}\Omega_{Y}\map{}\Omega_{X,G}$ and therefore a morphism
\begin{eqnarray*}
\inv{G}{d\pi} & : &\Omega_{Y}\map{}\inv{G}{\Omega_{X,G}}
\end{eqnarray*}
of coherent modules on $Y$. Under the additional assumption that $X$ is
smooth, then $\inv{G}{\Omega_{X,G}}$ is a torsion-free module and by (\cite[Theorem 1]{MR98k:14067}) the morphism
$\inv{G}{d\pi}$ is generically an isomorphism. Consequently, the kernel of
$\inv{G}{d\pi}$ is exactly the torsion of $\Omega_{Y}$ and we have an
inclusion~: $\bar\Omega_{Y}\subseteq\inv{G}{\Omega_{X,G}}$.
\end{remark}   

We now give some elementary properties of this construction~:
\begin{lemma}\label{pback} Let $f : X\map{}Y$ be an equivariant map of
$G$-schemes. Assume that the adjoint morphism
$\Omega_{Y}\map{}f_{*}\Omega_{X}$ is injective. Then the diagram~: 
$$
\xymatrix{
\Omega_{Y,G}\ar[r]\ar[d] & \Omega_{Y}\ar[d]\\
f_{*}\Omega_{X,G}\ar[r] & f_{*}\Omega_{X}
}
$$
is a fiber product diagram where all the morphisms are injective.
\end{lemma}
In other words, under the assumption, a differential form is horizontal 
if and only if its pull-back is.

\begin{proof}{\ref{pback}} The statement is an easy consequence of
the commutative diagram
$$
\xymatrix@+1.3em{
0\ar[r] & \Omega_{Y,G}\ar[r]\ar[d] & \Omega_{Y}\ar[r]^{d\mu_{Y,G}}\ar[d] &
 \Omega_{Y}\otimes\lieg^{\vee}\ar[d]  \\
0\ar[r] & f_{*}\Omega_{X,G}\ar[r] & f_{*}\Omega_{X}\ar[r]^{f_{*}d\mu_{X,G}} &
 f_{*}\Omega_{X}\otimes\lieg^{\vee} 
}
$$
where the two vertical morphisms on the left are injective by assumption.
\end{proof}

\begin{lemma}\label{fibration} Let $G$ be an algebraic group and $f :
X\map{}Y$ be a principal $G$-fibration.
Then the natural morphism 
$df : f^{*}\Omega_{Y}\map{}\Omega_{X,G}$ is an  isomorphism.
\end{lemma}
One is reduced to proving the statement in  the case of a trivial
$G$-fibration where this is obvious.

\begin{lemma}\label{twogroups} Let $G$ and $H$ be algebraic groups
acting on a scheme $X$. The natural commutative diagram
$$
\xymatrix @+1.3em{
0\ar[r] & \Omega_{X,G}\otimes\lieh^{\vee}\ar[r] &
 \Omega_{X}\otimes\lieh^{\vee}\ar[r]^{d\mu_{X,G}\otimes\lieh^{\vee}} & 
 \Omega_{X}\otimes\lieg^{\vee}\otimes\lieh^{\vee} \\
0\ar[r] & \Omega_{X,G}\ar[r]\ar[u] &
 \Omega_{X}\ar[r]^{d\mu_{X,G}}\ar[u]^{d\mu_{X,H}} & 
 \Omega_{X}\otimes\lieg^{\vee}\ar[u]^{-d\mu_{X,H}\otimes\lieg^{\vee}} \\
0\ar[r] & \Omega_{X,G}\cap\Omega_{X,H}\ar[r]\ar[u] &
 \Omega_{X,H}\ar[r]\ar[u] & 
 \Omega_{X,H}\otimes\lieg^{\vee}\ar[u] \\
 & 0\ar[u] & 0\ar[u] & 0\ar[u] 
}
$$
has exact rows and columns. Moreover, it induces an exact sequence~:
$$
\xymatrix @+1.3em{
0\ar[r] & \Omega_{X,G}\cap\Omega_{X,H}\ar[r] & \Omega_{X}\ar[r] & \Omega_{X}\otimes(\lieg^{\vee}\oplus\lieh^{\vee}).
}
$$
\end{lemma}
Observe, that we did not assume that the actions of $G$ and $H$ on $X$
commute, therefore this diagram is only separately $G$ and
$H$-equivariant, but, in general, not $G\times H$-equivariant.   

\section{The Euler derivation}
\label{euler}
We go on using the notations of section \ref{hordiff}.
Let $T=\Gm=\Spec(\field[\lambda,\lambda^{-1}])$ be a one-dimensional
torus with Lie algebra $\liet$ and $X$ an affine $T$-scheme.  
We recall that since $T$ is
abelian, the adjoint representation is trivial, i.e. $\liet$ is a
trivial $T$-module. We fix once for all an isomorphism
$\field\simeq\liet$ via the left-invariant
derivation~$\lambda\frac{\partial}{\partial\lambda}$. Composing the dual
of this last isomorphism with $d\mu^{1}_{X,T}$ we obtain a derivation
on $X$~:  
\begin{eqnarray*}
 \e_{X,T} & : & \Omega^{1}_{X}\map{}\O_{X}
\end{eqnarray*}
called the Euler derivation. 
Since $X$ is affine, we have $X=\Spec(A)$ with $A$ a graded ring. The
grading of $A$ corresponds to the weight for the $T$-action~: A section $f$
of $\O_{X}$ is said to be homogeneous of weight $w$ if
$\mu^{*}f=\lambda^{w}q^{*}f$. If $f$ is homogeneous of weight $w$, we set
$|f|=w$.

\begin{proposition}\label{euler1} Let $f$ be an homogeneous section of
$\O_{X}$. Then~:
$$
\e(df)=|f|f.
$$
\end{proposition}

\begin{proof}{\ref{euler1}} Let $w=|f|$. We have~:
\begin{eqnarray*}
\e(df) & = & \lambda\frac{\partial}{\partial\lambda}d\mu^{1}_{X,T}(df)\\
 & = &
\lambda\frac{\partial}{\partial\lambda}s^{*}(w\lambda^{w-1}f.d\lambda)\\
 & = &
\lambda\frac{\partial}{\partial\lambda}s^{*}(w\lambda^{w}f.\frac{d\lambda}{\lambda})\\ 
 & = & \lambda\frac{\partial}{\partial\lambda}(wf.\frac{d\lambda}{\lambda})\\
 & = & wf
\end{eqnarray*}
as expected.
\end{proof}

\begin{proposition}\label{eulerdef} The Euler derivation constructed above can be extended 
to a degree $-1$ endomorphism of the graded module $\Omega_{X}$ by 
setting~:
\begin{eqnarray*}
\e & : & \Omega^{k}_{X}\map{}\Omega^{k-1}_{X} \\
\e(df_{1}\wedge\ldots\wedge df_{k}) & = &
\sum_{i=1}^{k}(-1)^{k-i}\e(df_{i})df_{1}\wedge\ldots\wedge\widehat{df_{i}}\wedge\ldots\wedge df_{k}.  
\end{eqnarray*}
It satisfies the following two properties~:
\begin{enum}
    \item $\e^{2}=0$.
    \item For any two forms $\alpha,\beta$ of respective degree $k$ 
    and $l$, we have 
    $$\e(\alpha\wedge\beta)=(-1)^{l}\e(\alpha)\wedge\beta+\alpha\wedge\e(\beta).$$
\end{enum}
\end{proposition}
\begin{proof}{\ref{eulerdef}} By direct computation.
\end{proof}
We thus have constructed a complex that we will denote by $(\Omega_{X},\e)$.

The exterior differential algebra $(\Omega_{X},d)$ is also
graded by weight~: A section $\alpha$ of $(\Omega_{X},d)$ is
homogeneous of weight $w$ if $\mu^{*}\alpha=\lambda^{w}q^{*}\alpha$.
The following properties are then easy to check~:

\begin{proposition}\label{euler2} Let $\alpha$ and $\beta$ be homogeneous sections of
$\Omega_{X}$.
\begin{enum}
\item  The forms $d\alpha$ and $\e(\alpha)$ are
homogeneous and $|d\alpha|=|\e(\alpha)|=|\alpha|$.
\item The form $\alpha\wedge\beta$
is homogeneous and $|\alpha\wedge\beta|=|\alpha|+|\beta|$.
\item The algebra $\Omega_{X}$ is generated by the differentials of
homogeneous sections of $\O_{X}$.
\item $\Ker(\e)=\Omega_{X,T}$.
\end{enum}
\end{proposition}
  
\begin{proposition}\label{bracket} For any homogeneous $k$-forms $\alpha$, we have~:
$$
[\e,d]\alpha=(-1)^{k}|\alpha|\alpha.
$$
\end{proposition}
\begin{proof}{\ref{bracket}} This is a direct computation again.
\end{proof}
Let $c\geq 0$. The operators $\e$ and $d$ preserve the filtration by
codimension of the support and therefore they induce operators on
$\Gamma_{c}\Omega_{X}$ and $\Gamma_{(c)}\Omega_{X}$ that we again denote 
by $\e$ and $d$. Moreover, since  $\Gamma_{c}\Omega_{X}$ and
$\Gamma_{(c)}\Omega_{X}$ are also $T$-equivariant, the statement above
remains true for these modules.

\begin{proposition}\label{euler3} The submodule 
$\inv{T}{\Omega_{X,T}}\subseteq
\inv{T}{\Omega_{X}}$ is stable by the exterior
derivative of $\Omega_{X}$.
\end{proposition}
\begin{proof}{\ref{euler3}} Keeping in mind that $T$-invariants are
precisely homogeneous sections of null weight, the result is a direct
consequence of (\ref{bracket}) and (\ref{euler2}\,(iv)).  
\end{proof}

\section{Horizontal differentials~: Poincar\'e lemmas}

\begin{proposition}\label{dstable} Let $G$ be a reductive algebraic group
and $X$ an affine $G$-scheme.  Then the submodule $\inv{G}{\Omega_{X,G}}\subset
\inv{G}{\Omega_{X}}$ is stable by the exterior derivative of
$\Omega_{X}$. 
\end{proposition}

\begin{remark} This statement holds more generally for $G$ a linear
algebraic group. 
But its proof would require an algebraic construction of the Lie
derivative that we did not explain here. The proof would run as
follows~: For $v\in\lieg$, 
denotes by $\lie_{v}$ the Lie derivative and by $<v,\cdot>$ the
interior product. Then, for any section $\alpha$ of $\Omega_{X}$ we have
the relation~:
$$
\lie_{v}\alpha=d<v,\alpha>+<v,d\alpha>.
$$
The statement therefore follows from the observation that $\lie_{v}$
vanishes on $\inv{G}{\Omega_{X}}$.
\end{remark}
 
\begin{proof}{\ref{dstable}} We recall $\inv{G}{\Omega_{X}}$ is
obviously stable by exterior differentiation. Since $G$ is reductive, on can find
one-dimensional subtori $T_{1},\ldots,T_{d}$ of $G$ such that
$\lieg=\liet_{1}\oplus\ldots\oplus\liet_{d}$. Then, by
(\ref{twogroups}), we have~: 

\begin{eqnarray*}
\Omega_{X,G} & = &
 \Omega_{X,T_{1}}\cap\ldots\cap\Omega_{X,T_{d}}.
\end{eqnarray*}
And therefore

\begin{eqnarray*}
\inv{G}{\Omega_{X,G}} & = &
 \inv{G}{\Omega_{X}}\cap\Omega_{X,T_{1}}\cap\ldots\cap\Omega_{X,T_{d}}\\
 & = & \inv{G}{\Omega_{X}}\cap\inv{T_{1}}{\Omega_{X,T_{1}}}\cap\ldots\cap\inv{T_{d}}{\Omega_{X,T_{d}}}.
\end{eqnarray*}
By (\ref{euler3}), all the terms in the intersection above are
stable by $d$, so we can conclude that
$\inv{G}{\Omega_{X,G}}$ is stable by $d$ too.
\end{proof}

\begin{proposition}\label{estable} Let $G$ be a reductive algebraic
group and let $X$ be an affine $G\times T$-scheme. Then $\Omega_{X,G}$
is stable by $\e=\e_{X,T}$. We write $(\Omega_{X,G},\e)$ for this
subcomplex of $(\Omega_{X},\e)$.
\end{proposition}
\begin{proof}{\ref{estable}} By a direct calculation, using the 
explicit definitions of $d\mu_{X,G}$ and $\e$.
\end{proof}
Therefore, if $c\geq 0$ is an integer, $\Gamma_{c}\Omega_{X,G}$ is also
stable by $\e$ and therefore there is an induced endomorphism on
$\Gamma_{(c)}\Omega_{X,G}$.

\begin{corollary}\label{bracket2} Let $G$ be a reductive algebraic
group and let $X$ be an affine $G\times T$-scheme. Let $\alpha$ be a
homogeneous section (with respect to the $T$-action) of
$\inv{G}{\Omega^{k}_{X,G}}$. Then 
$$
[\e,d]\alpha=(-1)^{k}|\alpha|\alpha.
$$
\end{corollary}
Clearly, we again have a similar statement for $\Gamma_{c}\inv{G}{\Omega_{X,G}}$, $\Gamma_{(c)}\inv{G}{\Omega_{X,G}}$,
$\inv{G}{\Gamma_{c}\Omega_{X,G}}$ or $\inv{G}{\Gamma_{(c)}\Omega_{X,G}}$.

\begin{proposition}\label{euler5} Let $G$ be a reductive algebraic group and let $X$
be an affine $G\times T$-scheme. Then 
$$
\H{}{\inv{G}{\Omega_{X,G}},\e}=\H{}{\inv{G}{\Omega_{X,G}},\e}^{T}.   
$$
Let $c\geq 0$. Then the same relation holds for
$\Gamma_{c}\inv{G}{\Omega_{X,G}}$, $\Gamma_{(c)}\inv{G}{\Omega_{X,G}}$,
$\inv{G}{\Gamma_{c}\Omega_{X,G}}$ and $\inv{G}{\Gamma_{(c)}\Omega_{X,G}}$.
\end{proposition}

\begin{proof}{\ref{euler5}} Let $\alpha$ be a homogeneous section of
$\inv{G}{\Omega^{k}_{X,G}}\cap \Ker(\e)$. Then by (\ref{bracket2}) we
have $\e d\alpha=(-1)^{k}|\alpha|\alpha$. Therefore if $|\alpha|\neq 0$
the class of $\alpha$ in  $\H{k}{\inv{G}{\Omega_{X,G}},\e}$
vanishes. Since $\H{}{\inv{G}{\Omega_{X,G}},\e}^{T}$ is a direct factor
of $\H{}{\inv{G}{\Omega_{X,G}},\e}$, the equality is proved.
\end{proof}

\begin{lemma}\label{diffpur} Let $X$ be a quasi-conical affine
$T$-scheme. Then the pull-back morphism for the quotient map
$X\map{}X\quot T$ induces isomorphisms~:  
$$
\Omega_{X\quot T}\map{\sim}\inv{T}{\Omega_{X,T}}
\map{\sim}\inv{T}{\Omega_{X}}\subset\Omega_{X}.
$$
\end{lemma}

\begin{proof}{\ref{diffpur}} Easy, by arguments on  weights.
\end{proof}

\begin{proposition}\label{euler6} Let $G$ be a reductive algebraic
group and let $X$ be an affine $G\times T$-scheme, quasi-conical with
respect to the $T$-action. Then
the natural morphism
$$
\Omega_{X\quot T,G}\map{}\inv{T}{\Omega_{X,G}}
$$
induced by the $G$-equivariant map $X\map{}X\quot T$, is an
isomorphism.
\end{proposition}

\begin{proof}{\ref{euler6}} By (\ref{diffpur}) the hypotheses of
(\ref{pback}) are satisfied for the map $X\map{}X\quot T$. Taking
$T$-invariants in the diagram of (\ref{pback}) together with the
isomorphism $\Omega_{X\quot T}\map{\sim}\inv{T}{\Omega_{X}}$ 
gives the result.
\end{proof}

\begin{proposition}\label{euler8} Let $G$ be a reductive algebraic
group and let $X$ be an affine $G\times T$-scheme, quasi-conical with
respect to the $T$-action. Let
$d\geq c\geq 0$. There are isomorphisms of exact sequences
$$
\xymatrix{
0\ar[r] & \inv{T}{\Gamma_{d}\Omega_{X,G}}\ar[r]\ar@{=}[d] &
 \inv{T}{\Gamma_{c}\Omega_{X,G}}\ar[r]\ar@{=}[d] &
 \inv{T}{\Gamma_{c}/\Gamma_{d}\,\Omega_{X,G}}\ar[r]\ar@{=}[d] & 0\\
0\ar[r] & \H{}{\inv{T}{\Gamma_{d}\Omega_{X,G}},\e}\ar[r] &
 \H{}{\inv{T}{\Gamma_{c}\Omega_{X,G}},\e}\ar[r] &
 \H{}{\inv{T}{\Gamma_{c}/\Gamma_{d}\,\Omega_{X,G}},\e}\ar[r] & 0\\
0\ar[r] & \inv{G\times T}{\Gamma_{d}\Omega_{X,G}}\ar[r]\ar@{=}[d] &
 \inv{G\times T}{\Gamma_{c}\Omega_{X,G}}\ar[r]\ar@{=}[d] &
 \inv{G\times T}{\Gamma_{c}/\Gamma_{d}\,\Omega_{X,G}} \ar[r]\ar@{=}[d] & 0\\
0\ar[r] & \H{}{\inv{G}{\Gamma_{d}\Omega_{X,G}},\e}\ar[r] &
 \H{}{\inv{G}{\Gamma_{c}\Omega_{X,G}},\e}\ar[r] &
 \H{}{\inv{G}{\Gamma_{c}/\Gamma_{d}\,\Omega_{X,G}},\e}\ar[r] & 0
}
$$
\end{proposition}

\begin{proof}{\ref{euler8}} By (\ref{euler6}) we have
$\inv{T}{\Gamma_{c}\Omega_{X,G}}\subset\Omega_{X\quot T}$. Therefore
$\e$ vanishes for all the complexes involved in the first isomorphism
and this proves the first statement. 
For the second one, take $G$-invariants in the first diagram and use
(\ref{euler5}).
\end{proof}

One might understand the next two statements as a natural
generalisation, with $\e$ and $d$ exchanged, of the Poincar\'{e} Lemma
to singular varieties with reductive group action~:   

\begin{corollary}\label{euler9} Let $G$ be a reductive algebraic
group and let $X$ be an affine $G\times T$-scheme, quasi-conical with
respect to the $T$-action. Then
the $G$-equivariant map $X\map{}X\quot T$ induces an isomorphism
$$
\inv{G}{\Omega_{X\quot T,G}}\map{\sim}\H{}{\inv{G}{\Omega_{X,G}},\e}.
$$
\end{corollary}

\begin{corollary}\label{euler10} Let $G$ be a reductive algebraic
group and let $X$ be an affine $G\times T$-scheme, quasi-conical with
respect to the $T$-action and such that
$X\quot T=\Spec(\field)$. Then
$$
\H{}{\inv{G}{\Omega_{X,G}},\e}=\H{}{\inv{G}{\bar\Omega_{X,G}},\e}=\field.
$$
\end{corollary}
In particular, in the case of a trivial action
of $G$ on a variety and under the preceding hypotheses we have exact complexes
$$
\xymatrix{
\ldots\ar[r] &\Omega_{X}^{n}\ar[r] &\ldots\ar[r] &\Omega_{X}^{1}\ar[r]
&\O_{X}\ar[r] & \field\ar[r] & 0 \\
0\ar[r] &\bar{\Omega}_{X}^{n}\ar[r] &\ldots\ar[r] &\bar{\Omega}_{X}^{1}\ar[r]
&\O_{X}\ar[r] & \field\ar[r] & 0
}
$$

\section{Absolutely regular horizontal differentials}
\label{arhdiff}

In this section, we merge the construction of horizontal differentials
and the content of appendix \ref{aregdiff}.  

Let $X$ be a $G$-scheme and $f : \tilde X\map{} X$ a $G$-equivariant
desingularisation. We denote by $\tilde\Omega_{X,G}$ the sheaf
$f_{*}\Omega_{\tilde X,G}$. This definition is independent of the choice
of $f$, as in the non-equivariant case, since two equivariant
resolutions of singularities can be covered by a third one.

By construction, we have natural equivariant morphisms
$$
\Omega_{X,G}\map{}\tilde\Omega_{X,G}\map{}i_{*}\Omega_{X_{\smooth},G}
$$
where $i$ is the inclusion $X_{\smooth}\subset X$. Therefore, when $X$ is
reduced, we have~:
$$
\Omega_{X,G}\map{}\bar\Omega_{X,G}\subset\tilde\Omega_{X,G}\subset i_{*}\Omega_{X_{\smooth},G}.   
$$

\begin{proposition}\label{arhdiff2} Let $f : X\map{}Y$ be an equivariant dominant
morphism. Then we have a commutative diagram
$$
\xymatrix{
\Omega_{X,G}\ar[r] & \tilde\Omega_{X,G} \\
f^{*}\Omega_{Y,G}\ar[r]\ar[u] & f^{*}\tilde\Omega_{Y,G}\ar[u] 
}
$$
\end{proposition} 
 
\begin{proposition}\label{arhdiff3} Let $f : X\map{}Y$ be an invariant dominant
morphism. Then we have a commutative diagram
$$
\xymatrix{
\Omega_{X,G}\ar[r] & \tilde\Omega_{X,G} \\
f^{*}\Omega_{Y}\ar[r]\ar[u] & f^{*}\tilde\Omega_{Y}\ar[u] 
}
$$
\end{proposition}  

\begin{proposition}\label{arhdiff1} Let $f : X\map{}Y$ be a proper
equivariant birational morphism. Then the morphism
$\tilde\Omega_{Y,G}\map{}f_{*}\tilde\Omega_{X,G}$ is an isomorphism.  
\end{proposition}

With this at hand, we can give a partial answer to the question raised
by M. Brion (\cite[after Theorem 2]{MR98k:14067})~:

\begin{proposition}\label{arhcomp} Let $G$ be a reductive algebraic group, $X$ be
a smooth affine $G$-scheme and $\pi : X\map{}Y$ 
the quotient. There is a sequence of inclusions :
$$
\bar\Omega_{Y}\subseteq\tilde\Omega_{Y}\subseteq\inv{G}{\Omega_{X,G}}
\subseteq\omega_{Y} 
$$
which are equalities on the smooth locus of $Y$. 
\end{proposition}

\begin{corollary}\label{arhdual} Let $G$ be a reductive algebraic group, $X$ be
a smooth affine $G$-scheme with quotient map $\pi : X\map{}Y$ and let $n=\dim Y$. Then 
the dualizing sheaf on $Y$, $\omega^{n}_{Y}$, is isomorphic to 
$\inv{G}{\Omega^{n}_{X,G}}$.
\end{corollary}

\begin{proof}{\ref{arhcomp}} Since $\Omega_{X,G}=\tilde\Omega_{X,G}$, by
(\ref{arhdiff3}) we have inclusions
$\bar\Omega_{Y}\subseteq\tilde\Omega_{Y}\subseteq\inv{G}{\Omega_{X,G}}$
of torsion-free modules. Moreover, by the theorem of Brion
(\cite[Theorem 1]{MR98k:14067}), these are isomorphisms outside the closed subset
$Y_{\sing}$, therefore outside a closed subset of codimension $\geq
2$. Thus the modules involved have isomorphic biduals and we obtain~:
$$
\bar\Omega_{Y}\subseteq\tilde\Omega_{Y}\subseteq\inv{G}{\Omega_{X,G}}
\subseteq\dd{\Omega_{Y}}=\omega_{Y}. 
$$
\end{proof}

\begin{proof}{\ref{arhdual}} It is then a direct consequence of the fact that $Y$
has rational singularities (\cite{MR88a:14005}). Indeed, this implies
that $\tilde\Omega^{n}_{Y}\map{\sim}\omega^{n}_{Y}$. 
\end{proof}

\begin{remark}\label{arhdiff4} If one assume that all the points of $X$
are strongly 
stable for the action of $G$, i.e., that for all closed points $x\in X$,
the orbit $Gx$ is closed and the stabilizer $G_{x}$ is finite, then
there are isomorphisms
$$
\tilde{\Omega}_{Y}\map{\sim}\inv{G}{\Omega_{X,G}}
\map{\sim}\omega_{Y}.
$$

To prove this, one can assume that the group $G$ is already finite (use
the Etale Slice Theorem as in the last reduction step in
(\ref{ihdiffsmth}) below). With this assumption made it is easily seen
that $\Omega_{X,G}=\Omega_{X}$ (here $\lieg=(0)$) and that consequently
$\inv{G}{\Omega_{X}}=\omega_{Y}$. It therefore remains to see that
$\tilde{\Omega}_{Y}=\inv{G}{\Omega_{X}}$. This can be done as follows.

We have a commutative diagram
$$
\xymatrix{
\tilde X \ar[r]^{\tilde\pi}\ar[d]^{g} & \tilde Y\ar[d]^{f} \\
X\ar[r]^{\pi} & Y}
$$
where $f$ is a resolution of singularities for $Y$ and $\tilde X$ is the
normalization of the component birational to $X$ in $X\times_{Y}\tilde
Y$. The group $G$ acts naturally on $\tilde X$ and the map $\tilde\pi$
is the quotient morphism. We thus have a morphism 
$$
\Omega_{\tilde Y}\map{}\inv{G}{\tilde\pi_{*}\tilde{\Omega}_{\tilde X}}
$$
induced by $\tilde\pi$. Since $\tilde X$ is normal it is an isomorphism
in codimension 1 and since 
$\Omega_{\tilde Y}$ is locally free it is in fact an isomorphism (recalling
that $\tilde\Omega_{\tilde X}$ is torsion-free). Consequently, we have
$$
\tilde\Omega_{Y} = f_{*}\Omega_{\tilde Y}
 = f_{*}\inv{G}{\tilde\pi_{*}\tilde{\Omega}_{\tilde X}}
 = \inv{G}{\pi_{*}g_{*}\tilde{\Omega}_{\tilde X}} 
 = \inv{G}{\pi_{*}\Omega_{X}}.
$$
This proves our claim.
\end{remark}

\section{Invariant horizontal differentials and smoothness}
\label{ihdiffsmth}

In this section we give proofs for the results stated in the introduction~:

\begin{theorem}\label{smcrit1} Let $G$ be a reductive
algebraic group acting on a smooth affine variety $X$, with quotient 
map $\pi :X\map{} Y$ and let $k$ be an integer with $1\leq k\leq\dim Y$.
The morphism $\inv{G}{d\pi^{k}}$ is
surjective in codimension $k+1$ if and only if $Y$ is smooth in
codimension $k+1$.   
\end{theorem}

\begin{theorem}\label{smcrit2} Let $Y$ be the quotient of a smooth
affine variety by a reductive algebraic group and let $k$ be an integer with $1\leq k\leq\dim Y$. 
The fundamental class morphism 
$c_{Y}^{k}$ is surjective in codimension $k+1$
if and only if $Y$ is smooth in codimension $k+1$.   
\end{theorem}

\begin{proof}{\ref{smcrit1}} After deleting a closed subset of 
codimension $>k+1$ we may assume that the morphism
$\inv{G}{d\pi} : \Omega_{Y}\map{}\inv{G}{\Omega_{X,G}}$ is surjective in
degree $k$, i.e. that we have a surjection
$\Omega^{k}_{Y}\map{}\inv{G}{\Omega^{k}_{X,G}}$ and we want to prove
that under this hypothesis the singular locus of $Y$ has codimension $>k+1$. 

The proof, now divides in five steps.

\subsubsection*{Etale slices} 

Quite generally, let $H\map{}G$ be a map of
reductive algebraic groups and $W$ an affine $H$-scheme together with an
$H$-equivariant map $ j : W\map{} X$. We let 
$G\times H$ act on $G\times W$ in the following way : 
$(g,h)(g',w)=(gg'h^{-1},hw)$ and denote by $f : G\times W\map{}
G\times_{H}W$ the quotient by $1\times H$. Observe that since $1\times
H$ acts freely on $G\times W$, the map $f$ is a principal fibration and
therefore is smooth. We obtain commutative
diagram of $G\times H$-schemes~:

\begin{eqnarray}
\label{slicediag1}\xymatrix @+1.3em{
G\times W\ar[d]^{}\ar[r]^{f} &
 G\times_{H}W\ar[d]^{}\ar[r]^(.6){\bar{\mu}(G\times_{H} j)} & 
 X\ar[d]^{\pi}\\  
W\ar[r]^{} & W\quot H\ar[r]^{} & X\quot G
}
\end{eqnarray}
where the vertical maps are quotients by $G$, the horizontal maps in the
left-square are quotients by $1\times H$ and $\bar{\mu}$ is the
factorization of the $1\times H$-invariant map $\mu$ ($1\times H$ acts
trivially on $X)$.  

For $y\in Y$ a closed point, we denote by $T_{y}\subset X_{y}$ the
unique closed orbit over $y$. Let $x\in T_{y}$ be a closed point with (necessarily) reductive stabilizer
$H=G_{x}$. The Etale Slice theorem of Luna (\cite[pp 96--99]{MR49:7269}),
asserts the following~: There  
exists a smooth locally closed, $H$-stable subvariety $W$ of $X$ such
that $x\in W$, 
$G.W$ is an open set and such that in the natural commutative diagram
(\ref{slicediag1}) the right-square is cartesian with  etale  horizontal maps 
(i.e. an etale base change diagram). Moreover, letting $N=N_{T_{y}/X}(x)$ be the normal
space at $x$ of the orbit $T_{y}$, understood geometrically as a
rational representation of $H$, there is a natural map of $H$-schemes
$\rho : W\map{}N$, etale at $0$, which induces a commutative diagram~:

\begin{eqnarray}
\label{slicediag2}\xymatrix @+1.3em{
G\times_{H} N\ar[d]^{\phi} &
 G\times_{H}W\ar[d]^{}\ar[r]^(.6){\bar{\mu}(G\times_{H}
 j)}\ar[l]_{G\times_{H}\rho} & 
 X\ar[d]^{\pi}\\  
N\quot H & W\quot H\ar[r]^{}\ar[l] & X\quot G
}
\end{eqnarray}
where the two squares are cartesian and the horizontal maps are etale
neighbourhoods.

\subsubsection*{Stratification by slice type}
 
We again refer to (\cite[pp 100--102]{MR49:7269}). Let $H\subseteq G$ be a
reductive subgroup and $N$ an $H$-module. We have a commutative
diagram~:
\begin{eqnarray}\label{slicediag4}
\xymatrix @+1.3em{
G\times N\ar[r]^{f}\ar[rd] & G\times_{H} N\ar[d] \\
 & G/H
}
\end{eqnarray}
which realizes $G\times_{H}N$ as the total space of a $G$-equivariant
vector bundle over the affine homogeneous space $G/H$ with fiber at $1$
equals to $N$. Conversely let
$N$ be a $G$-equivariant vector bundle over an affine $G$-homogeneous
base $T$. Let $t\in T$ be a closed point then $N(t)$ is a $G_{t}$-module
and $G_{t}$ is reductive. Thus we have an equivalence between the set
$\{(H,N)\}$ up to conjugacy and the isomorphism classes of $G$-equivariant
vector bundles over affine homogeneous bases. We denote by 
${\mathcal M}(G)$ any of those sets and classes by brackets~$[\,]$. 

By the preceding, we thus have a map
$\mu : Y(\bar\field)\map{} {\mathcal M}(G)$ which sends $y$ to the
isomorphism class $[N_{T_{y}/X}\map{} T_{y}]$ or equivalently to the
``conjugacy class'' $[H, N]$ with the notations of the preceding section.
Let $\nu\in{\mathcal M}(G)$, then the set $\mu^{-1}(\nu)$ is a locally
closed subset of $Y$, smooth with its reduced scheme structure. We will
denote by $Y_{\nu}$ this smooth locally closed subscheme of
$Y$. Moreover the collection $\{Y_{\nu}\}_{\nu\in{\mathcal M}(G)}$ is a
finite stratification of $Y$ (in particular $\mu$ has finite
image). Therefore, the map $\mu$ can be extended to all the points of
$Y$~: Let $Z\subset Y$ be an irreducible closed subset, then there exists
a unique $\nu\in{\mathcal M}(G)$ such that $Z\cap Y_{\nu}$ is
dense in $Z$ and one can set $\mu(Z)=\nu$. Observe that  $\mu(Z)$ is the
slice type of a general point of $Z$. 

Another important fact about $\mu$ is that it is compatible with strongly
etale (also called excellent) morphisms~: Given such a map $\varphi$ between
smooth affine $G$-schemes, we have $\mu(\varphi\quot G)=\mu$. 

We now look closer to $G$-schemes of the kind $G\times_{H}N$ and their
quotients by $G$. Write $N_{H}$ for the canonical
complementary submodule to $N^{H}$ in $N$~: $N=N^{H}\times N_{H}$. Then
in the construction  of 
$G\times_{H}N$, $N^{H}$ is a trivial $H$-module and therefore the diagram obtained when $W$
is replaced by $N$ in the left square of (\ref{slicediag1})
reads~:
\begin{eqnarray}
\label{slicediag3}\xymatrix @+1.3em{
N^{H}\times(G\times N_{H}) \ar[d]^{p}\ar[r]^{f} &
 N^{H}\times(G\times_{H}N_{H})\ar[d]^{\phi} \\   
N^{H}\times N_{H}\ar[r]^{\psi} & N^{H}\times(N_{H}\quot H)
}
\end{eqnarray}
Let $\nu\in{\mathcal M}(G)$ be the class of $(H,N)$, then
$((G\times_{H}N)\quot G)_{\nu}=N^{H}\times 0\subseteq N\quot H$. One can
convince oneself of this fact through the description of $G\times_{H}N$
as an equivariant vector bundle over $G/H$.

\subsubsection*{Reduction to an isolated singularity} 

First, it is harmless to
assume that the singular locus of $Y$, $Y_{\sing}$ is irreducible.
Let $\mu(Y_{\sing})=\nu=[H,N]$ and let $y\in Y_{\sing}\cap Y_{\nu}$ be
a general closed point.
By standard
etale base change arguments in the diagram (\ref{slicediag2}), our hypothesis
and our conclusion hold for $\pi$ at $y$ if and only if they
respectively hold for
$\phi$ at  $0$. We can therefore assume that $X=G\times_{H}N$,
$\pi=\phi$ and $Y=N\quot H$. 

Now, with the notations of
(\ref{slicediag3}), it is clear that 
$Y_{\sing}=N^{H}\times(N_{H}\quot H)_{\sing}$. On the other hand 
$Y_{\nu}=N^{H}\times 0$ and, since
$\mu(Y_{\sing})=\nu$, the closed subset $Y_{\nu}$ should cut a dense open
set on $Y_{\sing}$. Consequently, we must have
$Y_{\nu}=Y_{\sing}$ and thus $(N_{H}\quot H)_{\sing}=0$. 

Let $\pi_{H} : X_{H}=G\times_{H}N_{H}\map{}Y_{H}=N_{H}\quot H$ be the
quotient map by $G$, then clearly $\pi=N^{H}\times \pi_{H}$. Let $k$ be
an integer, then the map $\inv{G}{d\pi}$ is diagonal with respect to the
decompositions~:
 
\begin{eqnarray*}
\inv{G}{\Omega^{k}_{X,G}} & = &
 \bigoplus_{i=0}^{k}\Omega^{i}_{N^{H}}\etimes\inv{G}{\Omega^{k-i}_{X_{H},G}}\\
\Omega^{k}_{Y} & = & 
 \bigoplus_{i=0}^{k}\Omega^{i}_{N^{H}}\etimes\,\Omega^{k-i}_{Y_{H}}
\end{eqnarray*}
Therefore $\inv{G}{d\pi}$ is surjective in degree $k$ if and only if
$\inv{G}{d\pi_{H}}$ is surjective in all degrees $k-\dim{N^{H}},\dots, k$.

To conclude, we can therefore make the extra assumption that 
$Y=X\quot G=N\quot H$ has only an isolated singularity at $0$. And one
should notice that 
the theorem remains in fact only to
be proved when $k=\dim Y-1$ or $\dim Y$, since, otherwise ($k<\dim
Y-1$) the statement is obviously true.

\subsubsection*{Reduction to the case of a representation}

We keep in mind all the identifications and assumptions made previously.
Recalling diagram (\ref{slicediag3}) and applying lemmas (\ref{fibration}) and (\ref{twogroups}) to the
fibration $f$, we have an exact sequence

\begin{eqnarray*}
\xymatrix @+1.3em{
0\ar[r] & f^{*}\Omega_{G\times_{H} N, G}\ar[r] & \Omega_{G\times N,
G}\ar[r] & \Omega_{G\times N, G}\otimes\lieh^{\vee}.
}
\end{eqnarray*}
Taking $G$-invariants together with lemma (\ref{fibration}) for $p$
leads to the exact sequence~:

\begin{eqnarray*}
\xymatrix @+1.3em{
0\ar[r] & \inv{G}{f^{*}\Omega_{G\times_{H} N, G}}\ar[r] &
\Omega_{N}\ar[r] & \Omega_{N}\otimes\lieh^{\vee} 
}
\end{eqnarray*}
Therefore, we have proved that 
$\inv{G}{f^{*}\Omega_{G\times_{H} N, G}}=\Omega_{N, H}$. Taking
$H$-invariants, we obtain 
$$
\inv{H}{\Omega_{N, H}}=\inv{G\times H}{f^{*}\Omega_{G\times_{H} N, G}}=
\inv{G}{\Omega_{G\times_{H} N, G}}.
$$
One can then conclude, that the hypothesis and the conclusion of the theorem
hold for $\phi$ if and only if they respectively hold for $\psi$. Thus
we are reduced to prove the theorem in the case where $X$ is a rational
representation of $G$ with $X\quot G$ having only an isolated
singularity at the origin.

\subsubsection*{Conclusion}

Carrying on, $X$ is now a rational $G$-module with quotient $\pi :
X\map{}Y$, such that $Y$ has only an isolated singularity at the origin. 
We recall the hypothesis in the theorem~: The morphism
$\inv{G}{d\pi}$ is surjective in degree $k\leq\dim Y$. We must prove
that $Y$ is smooth in codimension $k+1$. Thus we have to prove that if
$k=\dim Y$ or $\dim Y-1$ then $Y$ is smooth.

The one dimensional torus $T=\Gm$ acts on $X$ by homothety and this
action commutes with the action of $G$. Thus $X$ is a $G\times T$ scheme
and $Y$ is a $T$-scheme. Both $X$ and $Y$ are quasi-conical and 
$X\quot T=Y\quot T=\Spec(\field)$.

Let $n=\dim Y$. Applying (\ref{euler10}) to $X$ and $Y$ we obtain an
injective morphism of exact complexes (the kernel of $\inv{G}{d\pi}$ is
exactly the torsion of $\Omega_{Y}$, cf. remark \ref{torsionkernel})~:

$$
\xymatrix@-1em{
\inv{G}{\Omega_{X,G}} & 0\ar[r] &\inv{G}{\Omega_{X,G}^{n}}\ar[r]
 &\inv{G}{\Omega_{X,G}^{n-1}}\ar[r] &\ldots\ar[r]
 &\inv{G}{\Omega_{X,G}^{1}}\ar[r] 
 &\O_{Y}\ar[r] & \field\ar[r] & 0 \\
\bar{\Omega}_{Y}\ar[u]^{\inv{G}{d\pi}} &0\ar[r]
 &\bar{\Omega}_{Y}^{n}\ar[r]\ar[u] &\bar{\Omega}_{Y}^{n-1}\ar[r]\ar[u]
 &\ldots\ar[r] 
 &\bar{\Omega}_{Y}^{1}\ar[r]\ar[u]
 &\O_{Y}\ar[r]\ar@{=}[u] & \field\ar[r]\ar@{=}[u] & 0
}
$$
From this diagram, we deduce that if $\inv{G}{d\pi}$ is surjective in
degree $n-1$, then it is also surjective in degree $n$. Therefore we
have an isomorphism
$\bar\Omega^{n}_{Y}\map{\sim}\inv{G}{\Omega_{X,G}^{n}}$. Moreover, by
proposition (\ref{arhcomp}) we know that  
$\inv{G}{\Omega_{X,G}^{n}}=\omega^{n}_{Y}=\dd{\Omega^{n}_{Y}}$. Thus
$\bar\Omega^{n}_{Y}$ is a reflexive module. 

Recall that by the theorem of Boutot (\cite{MR88a:14005}), $Y$ has
rational singularities and in particular is normal and
Cohen-Macaulay and that $\omega^{n}_{Y}$ is then the dualizing module of
$Y$. The fundamental class map $c$ 
(\cite[5.2 p 91, 5.15 p 99]{MR90a:14021}, \cite{MR80h:14009} and
appendix \ref{rardiff}), in degree $n$, factors through~:

$$
\xymatrix{
\Omega^{n}_{Y}\ar[r]^{c}\ar@{>>}[d] & \omega^{n}_{Y} \\
\bar\Omega^{n}_{Y}\ar[ru] &  
}
$$
But $\bar\Omega^{n}_{Y}$ is reflexive and, since $Y$ is normal, $c$ is
an isomorphism in codimension $1$. Therefore $c$ is necessarily
surjective. 
We now invoke a theorem of Kunz and Waldi (\cite[5.22 p 107]{MR90a:14021})
to conclude that $Y$ is smooth.

The proof of theorem \ref{smcrit1} is complete.
\end{proof}

\begin{proof}{\ref{smcrit2}} Using (\ref{smcrit1}), we can a give a
straightforward proof of the result : By (\ref{arhcomp}) the hypotheses
of (\ref{smcrit1}) are satisfied for the same integer $k$.
\end{proof}

\section{The case of abelian finite groups}
\label{fgroup}

Let $G$ be a finite group acting on a quasi-projective scheme $X$ and
let $\pi : X\map{}Y$ be the quotient. 

For an element $g\in G$, we denote the closed subscheme of $g$-fixed points by
$X^{g}$ and for a point $x\in X$, we denote its stabilizer (also called
isotropy subgroup) by $G_{x}$. 
We then define a increasing filtration of $G$ by normal subgroups in the
following way~: For $k\geq 0$ an integer we set 
$G^{k}=<g\in G, \forall x\in X^{g}, \codim{X^{g},x}\leq k>$. In particular $G^{1}$ is the
subgroup generated by the {\em pseudo-reflections} in $G$. 
For a point $x\in X^{g}$, if $\codim{X^{g},x}\leq 1$ then $g$ is said to
be a {\em pseudo-reflection at $x$}. 
When $X$ is smooth, this condition is satisfied if and only
if locally at $x$, the diagonal form of $g$ is of the kind
$(\zeta,1,\ldots,1)$ for some root of unity $\zeta$. 
Clearly $g$ is a pseudo-reflection if and only if it is a
pseudo-reflection at all the points of $X^{g}$. 

When $G^{1}=(1)$ one says that $G$ is a {\em small} group of
automorphisms of $X$. In this case, by standard ramification theory, 
the quotient map is unramified in codimension one.
When $G=G^{1}$ one says that $G$ is generated by
pseudo-reflections. 
We now recall the classical

\begin{claim}{Theorem (Shephard-Todd, Chevalley, Serre \cite{MR15:600b,MR38:3267})} With the
preceding notations, the following conditions are equivalent :
\begin{enum}
\item The quotient $Y$ is smooth.
\item For all $x\in X$, the group $G_{x}$ is generated by
the pseudo-reflections at $x$.
\item The $\O_{Y}$-module $\pi_{*}\O_{X}$ is locally free.
\end{enum}
\end{claim}
Thus, the local study of  quotients of smooth varieties by finite groups
reduces to the study of quotients of smooth varieties by small finite
groups of automorphisms~: Indeed, the theorem above implies that, locally 
around $x$, the group 
$G_{x}/G_{x}^{1}$ is a small group of automorphisms of the smooth variety
$X/G_{x}^{1}$. 
It is also clear that, for local
questions, by the Etale Slice Theorem (see (\ref{ihdiffsmth}))
one is reduced to study the case where $X$ is a rational representation
of $G$. 

\begin{theorem}\label{fgroup3} Let $G$ be a finite abelian group 
acting on a smooth affine scheme $X$ with quotient $\pi : X\map{} Y$
and let $k$ be an integer with $1\leq k\leq\dim X$. The morphism
$\inv{G}{d\pi^{k}}$ is surjective if and only if $Y$ is smooth.
\end{theorem}

\begin{proof}{\ref{fgroup3}} By the preceding remarks, we are 
reduced to the case where $X$ is a rational representation of $G$ as 
a small group of automorphism. So
that the map $\pi$ is unramified in codimension one. 

We recall that, $G$ being finite, we have $\Omega_{X,G}=\Omega_{X}$. 
Moreover by (\ref{smcrit1}) we deduce that $Y$ is smooth in 
codimension 2 and we can assume  that $1\leq k <\dim{X}-1$. Thus we 
can assume that $\dim X>2$ and purity of 
the branch locus
implies that $\pi$ is unramified in codimension 2.

From now on we proceed by induction on $\dim X$. Since $G$ is 
abelian, $X$ decomposes as a product of representation~: $X=X'\times L$ 
with $2\leq\dim{X'}=\dim{X}-1$. We have a diagram
$$\xymatrix{
X'\ar@{^{(}->}[r] \ar[d]^{\pi'} & X\ar[d]^{\pi} \\
Y'\ar@{^{(}->}[r] & Y
}
$$
where the vertical maps are quotient by $G$ and the horizontal ones 
are embeddings. This induces a commutative diagram~:
$$\xymatrix{
\inv{G}{\Omega^{k}_{X}}\ar[r] & \inv{G}{\Omega^{k}_{X'}}\\
\Omega^{k}_{Y}\ar[r]\ar[u] & \Omega^{k}_{Y'}\ar[u]
}
$$
where all the morphisms are surjective. Thus, by the induction 
hypothesis, $Y'$ is smooth. Now, if $G$ were not trivial, the origin 
being a fixed point, the map $\pi'$ should have to be ramified
and, by purity of the branch locus again, its ramification locus 
should have codimension one.
But then $\pi$ should be ramified in codimension 2. It is a 
contradiction. Thus, $G$ is trivial and therefore $Y$ is smooth.
\end{proof}


\appendix

\section{Regular and absolutely regular differentials}
\label{rardiff}

\subsection{Regular differentials}
\label{regdiff}

Regular differentials together with duality theory have been studied by
many authors but from different viewpoints. The main results that we need
are found in the book of Kunz and Waldi (\cite{MR90a:14021}), but we
feel that the very general and explicit
construction of regular differentials in this book (where the
construction is local and relative from the beginning) asks a lot
of the (lazy) reader, and therefore does not ``specialize'' easily to a
convenient tool in the common case of schemes of finite type over a field. 

Thus we choose the following path~: We take the  theory of the residual
complex and fundamental class as exposed in the work of El Zein
(\cite{MR80h:14009}) as a ``black box'' and rephrase, with a view toward
Kunz and Waldi's theory of regular differentials, the  results and
constructions of El Zein. We do not intend to say anything new here and all the
subsequent claims are implicitely proved in El Zein's article
(\cite{MR80h:14009}). In fact, this approach was inspired to us by the
work of Kersken (\cite{MR85h:32014,MR86a:14015,MR85j:14032}). 

\subsubsection*{Construction}

Let $\field$ be a field of characteristic 0. For any scheme $X$ of finite
type over $\field$, there exists a {\em residual complex} $\K_{X}$
(\cite{MR36:5145}). This is a complex of injective
$\O_{X}$-modules concentrated in degree $[-\dim X, 0]$, the image of 
which in
the derived category is the {\em dualizing complex}.
 
Let $n=\dim X$. We denote by $\omega^{n}_{X}$ the module
$\cH{0}{\K_{X}[-n]}$. If $X$ is smooth, $\K_{X}$ is the Cousin resolution of
$\Omega^{n}_{X}[n]$. If $i : X\map{}Y$ is an embedding of $X$ into a
smooth $Y$ then
$\K_{X}=i^{!}\K_{Y}=\sHom_{\O_{Y}}(\O_{X},\K_{Y})$. If $\pi : X\map{}Y$
is a finite surjective morphism then the complexes $\K_{X}$ and
$\pi^{!}\K_{Y}$ are quasi-isomorphic and therefore $\omega^{n}_{X}\simeq\pi^{!}\omega^{n}_{Y}$. Moreover, the formation of the residual complex
commutes with restriction to an open set. 
Thus, for a general $X$, $\omega^{n}_{X}$ has the ${\mathrm S}_{2}$
property and coincides with $\Omega^{n}_{X}$ at the smooth
points of $X$. Consequently, if $X$ is normal then there is a natural
isomorphism  $\dd{\Omega^{n}_{X}}\map{\sim}\omega^{n}_{X}$.

The complex $\K_{X}$ is exact in degrees $\neq\dim X$ if and only if $X$
is equidimensionnal and Cohen-Macaulay. In this case, the module
$\omega^{n}_{X}$ is the {\em dualizing module} (usually denoted
$\omega_{X}$).

Now, following El Zein, let $\K_{X}^{*,\cdot}=\sHom(\Omega_{X},\K_{X})$. It
is a bigraded object, where the $*$ (resp.\ the $\cdot$) corresponds to
degrees in $\Omega_{X}$ (resp.\ in $\K_{X}$), concentrated in degrees
$[-\infty, 0]\times [-\dim X, 0]$. We now explain how one can put
on $\K_{X}^{*,\cdot}$ a structure of complex of right differential graded
$\Omega_{X}$-modules concentrated in degree $[-\dim X, 0]$.

The left $\Omega_{X}$-module structure of $\Omega_{X}$ given by exterior
product induces an obvious right $\Omega_{X}$-module structure on
$\K_{X}^{*,p}=\sHom(\Omega_{X},\K^{p}_{X})$ and the
differential $\delta$ of 
$\K_{X}$ induces an $\Omega_{X}$-linear differential~:
$\delta'=\sHom(\Omega_{X},\delta)$.  

The non-trivial point is the existence for all $p$ of a differential endo-operator $d'$
of order $\leq 1$ and $*$-degree 1 on
$\K_{X}^{*,p}$ satisfying the conditions
\begin{enum}
\item $\delta'.d'=d'.\delta'$.
\item $d'(\phi.\alpha)=\phi.(d\alpha)+(-1)^{q}(\d'\phi).\alpha$, for
$\alpha\in\Omega^{q}_{X}$ and $\phi\in\K_{X}^{*,p}$.
\end{enum} 
The construction of $d'$ is explained in (\cite[2.1.2]{MR80h:14009}),
the proof of (ii) follows from the lemma
(\cite[2.1.2, Lemme]{MR80h:14009}, be aware that there is a misprint in
this paper~: The logical section 2.1.2 is labelled 3.1.2) and 
the remarks following the proof of this lemma. Finally, (i) is a direct
consequence of  (\cite[2.1, Proposition]{MR80h:14009}) and 
(\cite[2.1.2, Proposition]{MR80h:14009}).  
We want to insist on the fact that, even in the smooth case, the operator $d'$ is not
the naive (and above all, meaningless) ``$\sHom(d,\K_{X})$''.
We can now define the module of {\em regular differential forms}~:
$\omega_{X}=\cH{*,0}{\K_{X}^{*,\cdot}[-n,-n]}$. Thus, $\omega_{X}$ is a right
differential graded $\Omega_{X}$-module and one has
$\omega^{k}_{X}=\sHom(\Omega^{n-k}_{X},\omega^{n}_{X})$.  

When $X$ is normal and equidimensional, the isomorphism
$\dd{\Omega^{n}_{X}}\map{\sim}\omega^{n}_{X}$ therefore induces an isomorphism
$\dd{\Omega_{X}}\map{\sim}\omega_{X}$. Thus, in this case, it is
easily seen that this construction coincides with that of Kunz and Waldi
(\cite[3.17, Theorem]{MR90a:14021}). Note also that, when $X$ is normal,
$\omega_{X}$ is a reflexive module. 

\subsubsection*{The fundamental class}

The {\em fundamental class} is constructed and studied by El Zein in
(\cite[3.1, Th\'{e}or\`{e}me]{MR80h:14009}). The fundamental class is
defined as a global section ${\mathrm C}_{X}$ of $\K^{*,\cdot}_{X}$ (as
a bigraded object) satisfying $d'{\mathrm C}_{X}=\delta'{\mathrm
C}_{X}=0$.
When $X$ is equidimensional of dimension $n$, the fundamental class is
homogeneous of degree $(-n, -n)$. In general, the contribution to
${\mathrm C}_{X}$ of an $m$-dimensional irreducible component of $X$ is
homogeneous of degree $(-m, -m)$ (cf. the next section). 
 Let $X$ be an $n$-dimensionnal scheme. By this observation, since
$\delta'{\mathrm C}_{X}=0$, we have an induced cohomology class
${\mathrm c}_{X}\in\omega^{0}_{X}$. 
Then, right multiplication defines a morphism
\begin{eqnarray*}
\Omega_{X} & \map{} & \omega_{X} \\
\alpha & \longmapsto & {\mathrm c}_{X}.\alpha 
\end{eqnarray*}
of differential graded $\Omega_{X}$-modules, thanks to
the relation $d'{\mathrm c}_{X}=0$. We again denote by ${\mathrm c}_{X}$ this
morphism and also call it the fundamental class morphism. 

To be a little more explicit,
${\mathrm c}_{X}\in\cH{0}{X, \K^{*,\cdot}_{X}[-n,-n]}=
\Hom(\Omega^{n}_{X}, \omega^{n}_{X})$ and the fundamental
class morphism in degree $k$ is the composition
$$
\Omega^{k}_{X}\map{}\Hom(\Omega^{n-k},\Omega^{n}_{X})\map{}
\Hom(\Omega^{n-k},\omega^{n}_{X})\simeq\omega^{k}_{X}.
$$

When  $X$ is normal and equidimensional, the morphism ${\mathrm c}_{X}$
can be identified with the natural morphism
$\Omega_{X}\map{}\dd{\Omega_{X}}\simeq\omega_{X}$.

We can now state the following fundamental theorem of Kunz and Waldi~:

\begin{claim}{Theorem (\cite[5.22, p107]{MR90a:14021})} Let $X$ be an
equidimensional Cohen-Macaulay reduced scheme of finite type over
$\field$ and let $n=\dim X$. Then the support of 
$\Coker ({\mathrm c}_{X})^{n}$ is precisely the singular locus of $X$.
\end{claim}

\subsubsection*{The trace map for regular differentials}

Let $f : X\map{}Y$ be a proper morphism, then the trace morphism 
$\Tr f : f_{*}\K^{*,\cdot}_{X}\map{}\K^{*,\cdot}_{Y}$ is obtained by the
composition of the natural morphism $\Omega_{Y}\map{}f_{*}\Omega_{X}$
with the trace morphism for residual complexes
$f_{*}\K_{X}\map{}\K_{Y}$. We thus have a well
defined trace morphism $\Tr f : f_{*}\omega_{X}\map{}\omega_{Y}$
vanishing if $\dim X \neq\dim Y$.

Assume that $f$ is birational, i.e., that there exists a dense open
subset $V\subset Y$ such that the induced morphism $f^{-1}(V)\map{}V$ be
an isomorphism. Then, by (\cite[3.1, Th\'{e}or\`{e}me]{MR80h:14009}) the
trace morphism $\Tr f : f_{*}\K^{*,\cdot}_{X}\map{}\K^{*,\cdot}_{Y}$
sends ${\mathrm C}_{X}$ to ${\mathrm C}_{Y}$. Consequently, under these
hypotheses we have a commutative diagram~:
$$
\xymatrix{
f_{*}\Omega_{X}\ar[r]^{{\mathrm c}_{X}} & f_{*}\omega_{X}\ar[d]^{\Tr f} \\ 
\Omega_{Y}\ar[r]^{{\mathrm c}_{Y}}\ar[u] & \omega_{Y}
}
$$

Let $X$ be a scheme and $X_{1},\ldots,X_{k}$ its irreducible components
with their reduced structure and inclusions 
$j_{i} : X_{i}\subset X$. Then by construction 
(\cite[p37]{MR80h:14009}) we have that 
${\mathrm C}_{X}=
\sum_{i} e_{X_{i}}(X)\Tr j_{i}({\mathrm C}_{X_{i}})$,
where $e_{X_{i}}(X)=\mathrm{length}(\O_{X,X_{i}})$, the multiplicity of
$X$ along $X_{i}$. 
Thus, we have     
${\mathrm c}_{X}=
\sum_{i} e_{X_{i}}(X)\Tr j_{i}({\mathrm c}_{X_{i}})$. 

Assume now that $f : X\map{}Y$ is a finite dominant morphism between
integral schemes then by (\cite[3.1, Proposition 2]{MR80h:14009}) we have that
$\Tr f({\mathrm C}_{X})=\deg(f){\mathrm C}_{Y}$. 
We therefore have a commutative diagram~: 
$$
\xymatrix{
f_{*}\Omega_{X}\ar[r]^{{\mathrm c}_{X}} & f_{*}\omega_{X}\ar[d]^{\Tr f} \\ 
\Omega_{Y}\ar[r]^{\deg(f){\mathrm c}_{Y}}\ar[u] & \omega_{Y} 
}
$$

\subsection{Absolutely regular differentials}
\label{aregdiff}

Let $X$ be a scheme and $f : \tilde X\map{}X$ a desingularisation 
(if $X$ is not reduced, by this, we mean a desingularisation of
$X_{\red}$). We recall that the $\O_{X}$-module $f_{*}\Omega_{\tilde X}$
is independent 
of the choice of $f$, we denote it by $\tilde\Omega_{X}$. It is usually
called the module of {\em absolutely regular differentials}, or
sometimes, when $X$ is a normal variety, the module of {\em Zariski
differentials}. By construction, we have natural morphisms
$$
\Omega_{X}\map{}\tilde\Omega_{X}\map{}i_{*}\Omega_{X_{\smooth}}
$$
where $i$ is the inclusion $X_{\smooth}\subset X$. Therefore, when $X$ is
reduced, we have~:
$$
\Omega_{X}\map{}\bar\Omega_{X}\subset\tilde\Omega_{X}\subset i_{*}\Omega_{X_{\smooth}}.   
$$
In general, we also have a commutative diagram~:
$$
\xymatrix{
f_{*}\Omega_{\tilde X}\ar@{=}[r] & f_{*}\omega_{\tilde X}\ar[d]^{\Tr f} \\
\Omega_{X}\ar[u]\ar[r]^{{\mathrm c}_{X}} & \omega_{X}
}
$$
and consequently, a sequence of morphisms
$$
\Omega_{X}\map{}\tilde\Omega_{X}\map{}\omega_{X}.
$$

Let $f : X\map{}Y$ be a dominant morphism. Then we
have a commutative diagram
$$
\xymatrix{
\Omega_{X}\ar[r] & \tilde\Omega_{X} \\
f^{*}\Omega_{Y}\ar[r]\ar[u] & f^{*}\tilde\Omega_{Y}\ar[u] 
}
$$
Assume moreover that the morphism $f$ is proper and birational. 
Then we have a commutative diagram
$$
\xymatrix{
f_{*}\Omega_{X}\ar[r] & f_{*}\tilde\Omega_{X}\ar[r] &
 f_{*}\omega_{X}\ar[d]^{\Tr f} \\ 
\Omega_{Y}\ar[r]\ar[u] & \tilde\Omega_{Y}\ar[u]\ar[r] & \omega_{Y}
}
$$
where the rows are factorisations of the respective fundamental class 
morphisms. Note that--obviously--the middle vertical arrow is an 
isomorphism.


\newpage
\section*{\refname}\thispagestyle{plain}

\bibliography{mybib}
\bibliographystyle{alpha}

\vfill

\hfill\vbox{
\hbox{Guillaume Jamet}\smallskip
\hbox{Universit\'{e} Pierre et Marie Curie}
\hbox{Institut de Math\'{e}matiques}
\hbox{Topologie et G\'eom\'etrie Alg\'ebrique, Case 247}
\hbox{4, place Jussieu}
\hbox{F-75252 PARIS Cedex 05}\smallskip
\hbox{T\'{e}l~: +33 1 44 27 86 58}
\hbox{e-mail~: {\tt jamet@math.jussieu.fr}}
}

\end{document}